\documentclass[11pt,reqno]{article}
\usepackage{amssymb, amsmath, amsthm, amsfonts, amscd, epsfig}
\usepackage[english]{babel}
\usepackage{graphicx, epsfig}
\usepackage{color}

\newtheorem{theorem}{Theorem}[section]

\newtheorem{lemma}[theorem]{Lemma}

\newtheorem{thm}{Theorem}[section]

\newtheorem{remark}{Remark}

\def\ep{\epsilon}

\newcommand{\nm}{\noalign{\smallskip}}
\newcommand{\ds}{\displaystyle}

\newcommand{\p}{\partial}
\newcommand{\pd}[2]{\frac {\p #1}{\p #2}}
\newcommand{\eqnref}[1]{(\ref {#1})}

\newcommand{\la}{\langle}
\newcommand{\ra}{\rangle}

\newcommand{\Kcal}{\mathcal{K}}

\newcommand{\Scal}{\mathcal{S}}

\newcommand{\beq}{\begin{equation}}
\newcommand{\eeq}{\end{equation}}

\numberwithin{equation}{section}
\numberwithin{figure}{section}

\begin{document}

\title{Field expansions for systems of strongly coupled plasmonic nanoparticles\thanks{\footnotesize The work of Hai Zhang was partially supported by HK RGC GRF grant 16304517  and 
startup fund R9355 from HKUST.}}

\author{Habib Ammari\thanks{\footnotesize Department of Mathematics, 
ETH Z\"urich, 
R\"amistrasse 101, CH-8092 Z\"urich, Switzerland (habib.ammari@math.ethz.ch, sanghyeon.yu@sam.math.ethz.ch).}  
\and  Matias Ruiz\thanks{\footnotesize Department of Mathematics and Applications,
Ecole Normale Sup\'erieure, 45 Rue d'Ulm, 75005 Paris, France
(matias.ruiz@ens.fr).} 
\and Sanghyeon Yu\footnotemark[2]   \and Hai Zhang\thanks{\footnotesize 
Department of Mathematics, 
HKUST, Clear Water Bay, Kowloon, Hong Kong (haizhang@ust.hk).}}

\date{} 

\maketitle

\begin{abstract}
This paper is concerned with efficient representations and approximations of the solution to the scattering problem by a system of strongly coupled plasmonic particles. Three schemes are developed: the first is the resonant expansion which uses the resonant modes of the system of particles computed by a conformal transformation, the second is the hybridized resonant expansion which uses linear combinations of the resonant modes for each of the particles in the system as a basis to represent the solution, and the last one is the multipole expansion with respect to the origin. By considering a system formed by two plasmonic particles of circular shape, we demonstrate the relations between these expansion schemes and their advantages and disadvantages both analytically and numerically. In particular, we emphasize the efficiency of the resonant expansion scheme in approximating the near field of the system of particles. The difference between these plasmonic particle systems and the nonresonant dielectric particle system is also highlighted. The paper provides a guidance on the challenges for numerical simulations of strongly coupled plasmonic systems. 

\end{abstract}




\section{Introduction}

Plasmonic particles are metallic particles with size in the range from several nanometers to hundreds of nanometers. Under the illumination of light, the free-electrons in the particle can be strongly coupled to light for certain frequencies resulting in strong scattering and enhancement of local fields. This phenomenon is called surface plasmon resonance \cite{SC10, Gri12} and the associated frequencies are called plasmonic resonant frequencies. Plasmonic particles has many applications in the field such as super-resolution in imaging \cite{simovski, matias}, super-focusing of light \cite{dmitri2013}, plasmonic bio-sensing \cite{anker, hh}. 
Because of plasmonic resonances, plasmonic particles are ideal subwavelength resonators for light and hence a basic building block for optical metamaterials and photonic crystals. 
The plasmonic resonance for a single particle has been intensively studied in the literature; see \cite{kelly, plasmon1, plasmon4} and the references therein from the physics and experimental aspect, and \cite{pierre, matias, matias2,kang1, Gri12} from mathematical aspect.

In this paper, we are concerned with the scattering by a system of strongly coupled plasmonic particles, where the separation distance between neighboring particles is comparable or smaller than the characteristic size of the particles. 
Compared to the one particle system, the multiple particle system offers a great opportunity to tune the plasmonic resonances to a wider frequency regime as well as more flexibility to control of propagation of light \cite{lauchner}. New phenomena such as Fano-resonance \cite{boris2010} and artificial optical magnetism \cite{alu} may occur and these lead to new applications. An efficient and accurate computation of the scattering (optical response) of such complex system is the key to the modeling and design in such applications.

The numerical solution of the scattering by multiple strongly interacting plasmonic particles is a challenging task. First, the scattered field is multiscale. More precisely, the field is much stronger in certain localized regions such as at the boundaries of the particles and in the small gap regions between neighboring particles than in the other regions. An accurate computation of such field distribution is particularly important in sensing applications. Second, the scattered field is very sensitive to the operating frequency. Due to the many resonances induced by the strong coupling of particles, the scattered field can be changed dramatically even over a small range of frequencies. In a recent work \cite{yu17}, by combining the technique of transformation optics and the image charge method, an asymptotic formula as well as an efficient numerical scheme was developed for a system of two closely touching spherical plasmonic particles in the three dimensional case.

We remark that the scattering properties of strongly coupled plasmonic particles are very different from those associated with a system of dielectric particles. In \cite{sinum}, it was shown that  
closely spaced dielectric particles can be approximated in the far field by an equivalent ellipse with an equivalent dielectric property having the same polarization tensor. Such a simple result does not hold in the case of plasmonic particles.  Moreover, the dipole approximation of closely spaced plasmonic particles is not enough to accurately approximate the scattered field. High-order multipoles are required.  

In this paper, three field expansion formulas are proposed to represent the scattered field, based on which numerical approximation schemes can be generated. The first one (formula (\ref{EF1}))  is the resonant expansion of the particle systems. In this scheme, the system of particles is treated as a single particle with a complicated shape of multiple connected components. 
We consider two nearly touching disks.
To reveal the complex nature of strong interaction between two disks, we use a conformal transformation. 
Then all the plasmonic resonances and resonant modes can be computed and the latter form a basis to represent the scattered field. 
The second one is the plsamon hybridization method for the multi-particle system (formula (\ref{EF3})) \cite{hybrid}. In this scheme, the resonances and resonant modes for each of the particles in the system are first computed. Then a basis which consists of linear combinations of these individual modes based on their interactions can be generated and used in the field expansion of the particle system. Beside these two schemes, we also provide the multipole expansion with respect to the origin (formula (\ref{EF2})) of the scattered field, which has been used for the scattering of localized inhomogenities in a background homogeneous media \cite{book3}. Using this expansion, we shall explain the reason why the field generated by two plasmonic disks cannot be approximated by that of an equivalent ellipse.  

Of particular  interest is the evaluation of the performances of the described field expansion schemes. Our paper is intended to highlight the challenges of numerical approximations for the strongly coupled plasmonic particle systems. By using an example of a system of two circular plasmonic particles, we demonstrate the relations between the proposed three expansion schemes and enumerate their advantages and disadvantages both analytically and numerically. An emphasis is placed on the efficiency of the resonant expansion in approximating the near field of the system of plasmonic particles.
 
We remark that there is a vast literature on the scattering problem from multiple plasmonic particle systems; see for instance \cite{hybrid, simovski}. We also note that the case of a system of weakly coupled plasmonic particles was considered in \cite{matias}. Based on perturbation theory, its analysis leads to a mathematical framework for understanding the superfocuing phenomena in resonant structures.

The paper is organized as follows. In section \ref{sec1} we formulate the plasmon resonance problem for a system of particles. Section \ref{sec2} is devoted to explicit computations for the case of two disks. In section \ref{sec_plasmon}, we consider the plasmon hybridization model where the plasmon modes of multiple particles are expressed in terms of interactions between the plasmon resonances of single particles. Using this model, the scattered field is represented as the sum of excited modes. 
In section \ref{sec_far},  we introduce the multipole expansion of the scattered field with respect to the origin. The expansion makes use of the concept of contracted generalized polarization tensors (CGPTs) \cite{book3}. Explicit calculations of the    CGPTs associated with two disks are derived. 
 In section \ref{sec6} we present a variety of numerical and compare the performances of the proposed numerical schemes.

\section{Plasmon resonances for multiple nanoparticles} \label{sec1}

\subsection{Problem formulation}

Let $B_1$ and $B_2$ be two bounded domains in $\mathbb{R}^2$ with smooth boundaries$\partial B_1$ and $\partial B_2$. We assume that $B_1$ and $B_2$ are disjoint. Let $\ep$ be the permittivity distribution given by
\beq\label{sigma:def}
\ep=\ep_m\chi(B_1)+\ep_m\chi(B_2)+\chi(\mathbb{R}^2\setminus (B_1\cup B_2),
\eeq
where $\chi$ denotes the characteristic function.

Given a harmonic function $H$ in the whole space $\mathbb{R}^2$, we consider the following transmission problem:
\beq
    \begin{cases}
        \nabla \cdot (\ep \nabla u) = 0 &\text{ in }\; \mathbb{R}^2, \\[1.5mm]
         u - H = O(|x|^{-1}) &\text{ as }\; |x| \rightarrow \infty.
    \end{cases}
    \label{transmission}
\eeq

The solution $u$ can be represented using boundary integral operators.

Define the Neumann-Poincar\'{e} operator $\mathcal{K}_{B_j}^*$ associated with $B_j, j=1,2$ by
$$
\mathcal{K}_{B_j}^* [\varphi](x) =\frac{1}{2\pi} \int_{\p B_j }  \frac{\langle x-y,\nu^{(j)}(x)\rangle}{|x-y|^2} \varphi(y) d\sigma(y) ,   \quad x \in \p B_j,
$$
and the single layer potential $\mathcal{S}_{B_j}$ by 
$$
\mathcal{S}_{B_j} [\varphi](x) =\frac{1}{2\pi} \int_{\p B_j }  \ln|x-y| \varphi(y) d\sigma(y) ,   \quad x \in \p B_j.
$$
Here, $\nu^{(j)}$ is the outward normal to $\partial B_j$. 

The Neumann-Poincar\'{e} type operator $\mathbb{K}^*$ associated with $B_1\cup B_2$ is  given by
$$\mathbb{K}^*:=
\left[
  \begin{array}{cc}
  \ds  \Kcal^*_{B_1} & \ds \pd{}{\nu^{(1)}}\Scal_{B_2} \\
  \ds \pd{}{\nu^{(2)}}\Scal_{B_1} & \ds\Kcal_{B_2}^* \\
   \end{array}
\right],$$

The solution $u$  admits the integral representation \cite{pierre, matias}
\beq
    u = H+ \Scal_{B_1}[\varphi_1] + \Scal_{B_2}[\varphi_2], 
    \label{scattered}
\eeq
where $(\varphi_1,\varphi_2)$ is the solution to the integral equation
$$
 (\lambda I - \mathbb{K}^*)\begin{bmatrix}\varphi_1\\  \varphi_2 \end{bmatrix} =\begin{bmatrix}\p_\nu H|_{\partial B_1}\\  \p_\nu H |_{\partial B_2} \end{bmatrix}
$$
and
\beq \label{deflambda} 
\lambda= \frac{\ep_m+1}{2(\ep_m-1)}.
\eeq  
The problem is to analyze the behavior of the solution $u$ when $\ep_m$ can be negative.

\subsection{Spectral decomposition of the Neumann-Poincar\'e operator}

Define the single layer potential $\mathbb{S}$ associated with $B_1\cup B_2$ by
$$\mathbb{S}=
\left[\begin{array}{cc}
\Scal_{B_1} |_{\p B_1}&\Scal_{B_2}|_{\p B_1}\\
\Scal_{B_1} |_{\p B_2}&\Scal_{B_2}|_{\p B_2}
\end{array}\right].$$
Let $H^{-1/2}_0(\p B_j)$ for $j=1,2$ be the zero mean subset of the standard Sobolev space 
$H^{-1/2}(\p B_j)$. 
Let $\mathcal{H}^*$ be $H^{-1/2}_0(\p B_1) \times H^{-1/2}_0(\p B_2)$ equipped with 
the inner product
\beq
\la \varphi,\psi\ra_*:=-\la\varphi,\mathbb{S}[\psi]\ra_{-1/2,1/2},\eeq
where $\la , \ra_{-1/2,1/2}$ denotes here the duality pairing between $H^{-1/2}(\p B_1) \times H^{-1/2}(\p B_2)$ and $H^{1/2}(\p B_1) \times H^{1/2}(\p B_2)$.
 
The Neumann-Poincar\'e  operator $\mathbb{K}^*$ is self-adjoint on $\mathcal{H}^*$. 
Since $\p B_1$ and $\p B_2$ are smooth,  $\mathbb{K}^*$ is compact. Therefore, $\mathbb{K}^*$ admits the spectral decomposition
\begin{equation} \label{sd}
\mathbb{K}^* = \sum_{n=1}^\infty \lambda_n \phi_n \otimes \phi_n,
\end{equation}
where $\otimes$ denotes the tensor product and $(\lambda_n, \phi_n)$ is the eigenvalue and the normalized eigenfunction pair of $\mathbb{K}^*$. 
Note that $|\lambda_n| \leq 1/2$; see \cite{ACKLM2}.

\subsection{Spectral representation of the solution $u$}
By the spectral decomposition (\ref{sd}), we have
$$
\varphi
=
\sum_{n=1}^\infty  \frac{ \langle \psi_H, \phi_n \rangle_*  }{\lambda- \lambda_n} \phi_n, 
$$
where
\beq\label{def_psi_H}
\psi_H = 
\begin{bmatrix}
\p_\nu H |_{\p B_1}
\\
 \p_\nu H |_{\p B_2}
\end{bmatrix}.
\eeq
For $\varphi = [\varphi_1, \varphi_2]^T$, where the superscript $T$ denotes the transpose, 
introduce
$$
\widetilde{\mathbb{S}}[\varphi](x) 
=\mathcal{S}_{B_1}[\varphi_1](x) +  \mathcal{S}_{B_2}[\varphi_2](x), \quad x\in \mathbb{R}^2. 
$$
Then we have
\begin{equation} \label{EF1n}
u-H = \sum_{n=1}^\infty  \frac{ \langle \psi_H, \phi_n \rangle_*  }{\lambda- \lambda_n} \widetilde{\mathbb{S}}[\phi_n].
\end{equation}

If the permittivity $\ep_m$ is negative, then $\lambda$ can be close to one of the eigenvalues $\lambda_n$ of $\mathbb{K}^*$. If $\lambda$ is close to $\lambda_n$, then the eigenfunction $\phi_n$ is amplified. Typically, the field distribution of $\widetilde{\mathbb{S}}[\phi_n]$ shows oscillating behavior and its gradient (the electric field) is large near the boundary of $\p B_1 \cup \p B_2$.  This explains the plasmonic resonance mechanism. 

The permittivity $\ep_m$ of the plasmonic particles is frequency-dependent and is modeled by
\beq\label{drude}
\ep_m= \ep_m(\omega) = 1-\frac{\omega_p^2}{\omega(\omega+ i\gamma )}.
\eeq
Here $\omega_p$ is the plasma frequency and $\gamma$ is the damping parameter. When $\omega<\omega_p$, the real part of $\ep_m$ is clearly negative. When $\gamma=0$, the frequency $\omega_n$ corresponding to eigenvalue $\lambda_n$ is called the resonant frequency.

\section{Conformal transformation based description} \label{sec2}

Here we give explicit computations of plasmon resonances for the case of two disks. We remark that the results in this section were obtained in \cite{mihai}. See also \cite{BT_disks}. Suppose that $B_1$ and $B_2$ are two disks with the same radius $R$ and let $d$ be the separation distance
$$d:=\mbox{dist}(B_1,B_2).$$
 We set the Cartesian coordinates $(x_1,x_2)$ such that the $x_1$-axis is parallel to the line joining the centers of the two disks. In other words, the disk $B_j$ is centered at $(-1)^i(R+d/2,0)$.

\subsection{Mobius transformation and the bipolar coordinates}

Each point $x=(x_1,x_2)$ in the Cartesian coordinate system corresponds to  $(\zeta,\eta)\in\mathbb{R}\times (-\pi,\pi]$ in the bipolar coordinate system through the equations
\begin{equation} \label{bipolar}
x_1=\alpha\frac{ \sinh \zeta }{\cosh \zeta - \cos \eta} \quad \mbox{ and } \quad x_2=\alpha\frac{\sin \eta}{\cosh \zeta - \cos \eta}
\end{equation}
with a positive number $\alpha$. Notice that the bipolar coordinates can be defined using a conformal mapping. Define the conformal map $\Psi$ by
$$
z=x_1 + i x_2= \Psi(\tilde{z})=\alpha \frac{\tilde z+1}{\tilde z-1}.
$$
If we write $\tilde z=e^{\zeta- i\eta}$, then we can recover \eqnref{bipolar}. 

From (\ref{bipolar}), we can derive that the coordinate curves  $\{\zeta=c\}$ and $\{\eta=c\}$ are, respectively, the zero-level set of the following two functions:
\beq\label{function:f}f_\zeta(x_1,x_2)=\left(x_1-\alpha\frac{\cosh c}{\sinh c}\right)^2 +x_2^2-\left(\frac{\alpha}{\sinh c}\right)^2\eeq
and $$f_{\eta}(x_1,x_2)=x_1^2 +\left(x_2-\alpha\frac{\cos c}{\sin c}\right)^2-\left(\frac{\alpha}{\sin c}\right)^2.$$\\

Let $\{\hat{\mathbf{e}}_{\zeta},\hat{\mathbf{e}}_{\eta}\}$ be the orthonormal basis vectors for the bipolar coordinates given by
$$
\hat{\mathbf{e}}_{\zeta}= \frac{\p \mathbf{x}/ \p \zeta}{|\p \mathbf{x}/ \p\zeta|} \quad \mbox{and} \quad \hat{\mathbf{e}}_{\eta}= \frac{\p \mathbf{x}/ \p \eta}{|\p \mathbf{x}/ \p\eta|}.
$$
We also denote the standard unit basis vectors in $\mathbb{R}^2$ by $\{\mathbf{e}_1,\mathbf{e}_2\}$.

In the bipolar coordinates, the scaling factor $h$ is 
$$
h(\zeta,\eta)=\frac{\cosh\zeta-\cos\eta}{\alpha}.
$$
The gradient of any scalar function $g$ is given by
\beq \label{grad_bipolar}
\nabla g = h(\zeta,\eta)\left( \frac{\p g}{\p\zeta}\hat{\mathbf{e}}_{\zeta}+ \frac{\p g}{\p\eta}\hat{\mathbf{e}}_{\eta}\right).
\eeq
Moreover, the normal and tangential derivatives of a function $u$ in bipolar coordinates are 
\begin{eqnarray}\label{nor_bipolar}
\left\{ \begin{array}{l} 
\ds \pd{u}{\nu}\Bigr|_{\zeta=c}=\nabla u\cdot \nu_{\zeta=c}=-\mbox{sgn}(c)h(c,\eta)\pd{u}{\zeta}\Bigr|_{\zeta=c},\\
\nm
 \ds \pd{u}{T}\Bigr|_{\zeta=c}=-\mbox{sgn}(c)h(c,\eta)\pd{u}{\eta}\Bigr|_{\zeta=c},
\end{array} \right. \end{eqnarray}
and the line element $d\sigma$ on the boundary $\{\zeta=s\}$ is
$$
d\sigma = \frac{1}{h(s,\eta)} d\eta.\\
$$
Here, $\mbox{sgn}$ denotes the sign function.

Using \eqnref{bipolar}, we have the following harmonic expansions for the two linear functions $x_1$ and $x_2$:
\beq\label{x1_bipolar}
x_1=\mbox{sgn}({\zeta})\alpha\left[1+2\sum_{n=1}^\infty  e^{-n|\zeta|}\cos n\eta\right],
\eeq
and
$$x_2=2\alpha\sum_{n=1}^\infty  e^{-n|\zeta|}\sin n\eta.$$

\subsection{Plasmon resonance modes for two separated disks}\label{NP_separated}

We consider the spectral properties of the Neumann-Poincar\'e operator $\mathbb{K}^*$ for two disks $B_1 \cup B_2$. 
Set
\beq\label{def_alpha_xi0}
\alpha= \sqrt{ d (R+ \frac{d}{4})}\quad\mbox{and}
\quad s=\sinh^{-1}\left(\frac{\alpha}{R}\right).
\eeq
It is easy to check that
\beq
\p B_j=\{\zeta=(-1)^js\},\quad \mbox{for }j=1,2.
\eeq
Note that the center of $B_j$ can be rewritten as $((-1)^{j}\alpha \coth s, 0)$.

The spectral properties can be represented in terms of the parameter $\alpha, s$ and the bipolar coordinates.
The following spectral decomposition of $\mathbb{K}^*$ was derived in \cite{mihai}.  
\begin{thm}\label{spectral_Kbbstar}
 We have the following spectral decomposition of $\mathbb{K}^*$ on $\mathcal{H}^*$:
\beq \label{sd2}
\mathbb{K}^* = \sum_{n\neq 0} \lambda_n^+ \Psi_{n}^+ \otimes \Psi_{n}^+
+  \lambda_n^- \Psi_{n}^- \otimes \Psi_{n}^-,
\eeq
where $\lambda^\pm_{n}$ are the eigenvalues of $\mathbb{K}^*$ given by
$$
\lambda^\pm_{n} = \pm\frac{ 1}{2}e^{-2|n|s}, \quad n\neq 0,
$$
and $\Psi_{n}^\pm$ are their associated (normalized) eigenfunctions defined by
\beq\label{normal_eigenvec}
\Psi_{n}^\pm (\eta)=  c_n^\pm {h(s,\eta) e^{in\eta}}
\begin{bmatrix}
1 \\
\mp 1
\end{bmatrix}, \quad c_n^\pm = \frac{\sqrt{|n|}}{\sqrt{4\pi (1/2-\lambda_n^\pm)}}, \quad n\neq 0.
\eeq

\end{thm}
The single layer potentials $\widetilde{\mathbb{S}}[\Psi_n^\pm]$ of the eigenfunctions can be explicitly computed as the following lemma.

\begin{lemma} We have the explicit formulas for the single layer potentials  $\widetilde{\mathbb{S}}[\Psi_n^\pm]$ of the eigenfunctions as follows:
\begin{align}
&\big( \widetilde{\mathbb{S}}[\Psi_n^\pm]\big)(\zeta,\eta)=(\mbox{const.})+
  \nonumber
 \\
 & \qquad \qquad \quad c_n^\pm\times
\begin{cases}
\ds \mp \frac{1}{2|n|}(e^{|n|s}\mp e^{-|n|s}) e^{|n|\zeta+in\eta},\quad&\mbox{for }\zeta<-s \  (\mbox{or } B_1),\\[3mm]
\ds \frac{-1}{2|n|}e^{-|n|s}(e^{-|n|\zeta}\mp e^{|n|\zeta})e^{in\eta},\quad&\mbox{for }-s<\zeta<s \  (\mbox{or } \mathbb{R}^2\setminus(B_1\cup B_2)),\\[3mm]
\ds \frac{-1}{2|n|}(e^{-|n|s}\mp e^{|n|s}) e^{-|n|\zeta+in\eta},\quad&\mbox{for }\zeta>s \  (\mbox{or } B_2).
\end{cases}
\label{Single_normal_eigen}
\end{align}
\end{lemma}

\subsection{Exact analytic solution for the scattered field using the spectral decomposition (\ref{sd2})}\label{subsec:exact}

Here we explicitly compute the solution $u$ when $H(x_1,x_2)=x_1$ using the spectral decomposition of the NP operator $\mathbb{K}^*$.

By the representation \eqnref{EF1n} and the spectral properties given in Theorem \ref{spectral_Kbbstar}, we have
\begin{align}
(u-H)(x) 
&=\sum_{n\neq 0} \frac{ \langle\psi_H, \Psi_n^+ \rangle_*  }{\lambda- \lambda_k^+} \widetilde{\mathbb{S}}[\Psi_n^+](x)
+
\sum_{n\neq 0} \frac{ \langle\psi_H, \Psi_n^- \rangle_*  }{\lambda- \lambda_k^-} \widetilde{\mathbb{S}}[\Psi_n^-](x)\nonumber
\\
&=
\sum_{n\neq 0} \frac{ (\frac{1}{2}-\lambda_n^+)\langle H, \Psi_n^+ \rangle_{L^2}  }{\lambda- \lambda_k^+} \widetilde{\mathbb{S}}[\Psi_n^+](x)
+
\sum_{n\neq 0} \frac{ (\frac{1}{2}-\lambda_n^-)\langle H, \Psi_n^- \rangle_{L^2}  }{\lambda- \lambda_k^-} \widetilde{\mathbb{S}}[\Psi_n^-](x),
\end{align}
where we have used
\begin{align*}
 \langle\psi_H, \Psi_n^\pm \rangle_* &=- \langle \p_\nu H, \mathbb{S}[\Psi_n^\pm] \rangle_{ L^2}=-\langle H, \p_\nu\mathbb{S}[\Psi_n^\pm] \rangle_{ L^2} 
 \\
 &= \langle H, (\frac{1}{2}-\mathbb{K}^*)[\Psi_n^+] \rangle_{L^2}
 =(\frac{1}{2}-\lambda_n^\pm) \langle H, \Psi_n^+ \rangle_{L^2}. 
\end{align*}

Suppose that $H(x_1,x_2)=x_1$. Since we have from \eqnref{x1_bipolar} that
$$
 H |_{\p D_i} = (-1)^{i}\alpha\Big(1+2\sum_{n=1}^\infty e^{-n s} \cos n \eta\Big),
$$
we can easily check that
$$
(\frac{1}{2}-\lambda_n^+)\langle H, \Psi_n^+ \rangle_{L^2} =  \left(\frac{|n|}{4\pi(c_n^+)^2} \right) (-4\pi \alpha  c_n^+ e^{-|n| s}) = (-c_n^+)^{-1} \alpha |n| e^{-|n|s},
$$
and
$$
\langle H, \Psi_n^- \rangle_{L^2} = 0.
$$
Thus, we get the following result.
\begin{theorem} When $H(x_1,x_2)=x_1$, the scattered field $u-H$ can be represented in terms of the eigenfunctions $\Psi_n^+$ (or the plasmon resonance modes for two disks $B_1\cup B_2$) as follows:
\beq\label{u_H_spectral}
(u-H)(x) = \sum_{n\neq 0} (-c_n^+)^{-1}\frac{ \alpha |n| e^{-|n|s} }{\lambda- \lambda_n^+} \widetilde{\mathbb{S}}[\Psi_n^+](x),
\eeq
where $c_n^+$ is given as \eqnref{normal_eigenvec}.
\end{theorem}

\begin{remark}
In fact, we can simplify (\ref{u_H_spectral}) in a more explicit form. By \eqnref{Single_normal_eigen}, we have
\begin{equation} \label{EF1}
(u-H)(\zeta,\eta) = \sum_{n=1}^\infty \frac{(-2)\alpha e^{-2n s}}{\lambda - \lambda_n^+} \sinh n\zeta \cos {n\eta}
\end{equation}
for $-s\leq \zeta \leq s \mbox{ (or } \mathbb{R}^2\setminus(B_1\cup B_2))$.
\end{remark}

\section{Plasmon hybridization description} \label{sec_plasmon}

 In the plasmon hybridization model \cite{hybrid},
 the plasmon resonance modes of multiple particles can be understood as a linear combination of the resonance modes of individual particles.
 In this section, we rigorously justify the plasmon hybridization model.
  This model has been used for the (heuristic) analysis of a variety of plasmonic systems.

\subsection{ Field expansion using plasmon resonances of a single disk}

We first consider the plasmon resonances of individual disk.
For the disk $B_j$, $0$ is the only eigenvalue of the Neumann-Poincar\'e  operator $\mathcal{K}_{B_j}^*$. Its associated eigenfunctions on $\mathcal{H}^*(\p B_j)$ are
$\phi_m^{(j)}:=e^{i m\theta_j}$ for $m \neq 0$. The single layer potentials generated by these eigenfunctions can be easily computed as
\begin{align}
\mathcal{S}_{B_j}\big[\phi_{m}^{(j)}\big] (r_j,\theta_j) = -\frac{R^{m+1}}{2m}  \frac{ e^{im \theta_j}}{r_j^m} , \quad r_j > R.
\label{single_layer_disk_exp}
\end{align}
These single layer potentials form a basis for the scattered field $u-H$ outside the two disks $B_1 \cup B_2$. 

When $H(x_1,x_2)=x_1$, the field $u-H$ can be represented as the following multipole expansion with respect to the center of each of the two disks.  
In fact, :
\begin{align}
(u-H)(x) = \sum_{m\neq 0} M^{(1)}_m \frac{e^{i m \theta_1}}{r_1^{m} } + M_m^{(2)} \frac{e^{i m \theta_2}}{r_2^{m} },\quad x\in \mathbb{R}^2\setminus(B_1\cup B_2),  \label{EF3}
\end{align}
with multipole coefficients $M^{(j)}_m$.
Note that $M_{ 1}^{(j)}$ (or $M_{ 2}^{(j)}$) is the dipole moment (or the quadrupole moment) induced on the disk $B_j, j=1,2$.

\subsection{Plasmon hybridization}

Consider the exact eigenfunction $\Psi_n^\pm $ for the two disks defined by \eqnref{normal_eigenvec}. We show that $\Psi_n^\pm $ can be decomposed in terms of $\phi_m^{(j)}=e^{im \theta_j}$, which are eigenfunctions for a single disk $B_j$, as follows:
$$
\Psi_n^\pm = \sum_{m\neq 0}  \begin{bmatrix}
 a_{nm}^{(1)} e^{i m\theta_1} \\ a_{nm}^{(2)} e^{i m\theta_2}
 \end{bmatrix}
$$
with some coefficients $a_{nm}^{(j)}$. We will compute $a_{nm}^{(j)}$ explicitly.

 First we represent $h(s,\eta)e^{i n \eta}$ as a Fourier series on $\p B_j$. Due to the symmetry of $B_1 \cup B_2$, we have 
\beq\label{TO_multipole1}
h(s,\eta)e^{in \eta} = \sum_{m\neq 0} b_{nm} e^{i m \theta_1} \quad \mbox{on }\p B_1,
\eeq
and
\beq\label{TO_multipole2}
h(s,\eta)e^{i n \eta} = \sum_{m\neq 0} (-1)^{m} b_{nm} e^{i m \theta_2} \quad \mbox{on }\p B_2
\eeq
with for some coefficients $b_{nm}$ satisfying $b_{nm}=b_{n(-m)}$.

We now compute $b_{nm}$ explicitly.
Recall that the center of $B_1$ is $(-\alpha \coth s, 0)$. Then, by the definition of the bipolar coordinates, we have for $z\in \p B_1$ (or $\zeta=-s$) that
\begin{align}
\frac{e^{-im \theta_1}}{R^{m}} = (z+\alpha \coth s)^{-m} = \left(\alpha \frac{e^{-s-i\eta}+1}{e^{-s-i\eta}-1}+\alpha \coth s\right)^{-m}.
\end{align}
Then, by straightforward but tedious computations, we obtain 
\begin{align}
\frac{e^{-i m \theta_1}}{R^{m}} = \frac{1}{\alpha^m(\coth s+1)^m }
\sum_{k=0}^\infty {F}(k,m,s)
e^{-ks} e^{-i k\eta}, \quad m \geq 1,
\end{align}
where
$$
{F}(k,m,s)=
\sum_{l=0}
^{\min(k,m)} 
(-1)^l
\begin{pmatrix}
m
\\
l
\end{pmatrix}
\begin{pmatrix}
m+k-l-1
\\
m-1
\end{pmatrix}
e^{2ls}.
$$
Then, since $b_{nm} = (2\pi R)^{-1}\int_{\p B_1} h(s,\eta)e^{i n\eta} e^{-i m \theta_1} d\sigma$ and $d\sigma = h(s,\eta)^{-1} d\eta$, we have
\beq\label{def_bnm}
b_{nm} = \frac{R^{|m|-1} e^{-ns}}{ 2\alpha^{|m|}(\coth s+1)^{|m|}}{F}(n,|m|,s).
\eeq

So we get the following result on the plasmon hybridyzation.

\begin{lemma} The eigenfunctions $\Psi_n^\pm$ for two separated disks $B_1\cup B_2$ can be represented as a linear combination of the eigenfunctions $e^{i m \theta_j}$ for a single disk $B_j$ as follows:
$$
\Psi_n^\pm =  \sum_{m\neq 0}  c_m^\pm b_{nm}  \begin{bmatrix}
 e^{i m\theta_1} 
 \\ \mp  (-1)^{m} e^{i m\theta_2}
 \end{bmatrix} 
$$
where $c_n^\pm$  and  $b_{nm}$ are given by \eqnref{normal_eigenvec} and  \eqnref{def_bnm}, respectively.
\end{lemma}

\begin{remark}
By assuming the distance $d$ is large and deriving the asymptotics of $a_{nm}$, we can easily get
\begin{equation} \label{HOT}
\Psi_n^\pm = c_n^+ b_{nn}\begin{bmatrix}
 \cos n\theta_1 \\ \mp (-1)^n  \cos n \theta_2
 \end{bmatrix} + O(d^{-1}).
\end{equation}
On the other hand, the eigenvalues are approximately given by 
\begin{equation} \label{HOT2}
\lambda_n^\pm \sim  \pm \frac{1}{2} d^{-2n}.
\end{equation}
Here, $+$ (resp. $-$) sign case  is called the bonding mode (resp. anti bonding mode). 
This plasmon hybridization model provides a simple description when the disks are well-separated. But, when the disks are nearly touching, the higher order terms in (\ref{HOT}) cannot be anymore neglected. Hence, the plasmon hybridization picture becomes very complicated in this case.
\end{remark}

\subsection{Computation of the multipole coefficients $M_{m}^{(j)}$}
Here we compute $M_m^{(j)}$ explicitly to describe the scattered field $u-H$ (when $H(x_1,x_2)=x_1$) using plasmon hybridization.

We have from 
From \eqnref{u_H_spectral} and \eqnref{single_layer_disk_exp} that
\begin{align}
(u-H)(x) &= \sum_{n\neq 0} (-c_n^+)^{-1}\frac{ \alpha |n| e^{-|n|s} }{\lambda- \lambda_n^+} \widetilde{\mathbb{S}}[\Psi_n^+](x)\nonumber
\\
&= \sum_{m\neq 0}\left(\sum_{n\neq 0}\frac{ (-1)\alpha |n| e^{-|n|s} }{\lambda- \lambda_n^+} b_{nm}\right)\big(\mathcal{S}_{B_1}[\phi_{m}^{(1)}] +(-1)^{m+1} \mathcal{S}_{B_2}[\phi_{m}^{(2)}]\big).\nonumber
\end{align}
Then, using \eqnref{single_layer_disk_exp} and \eqnref{EF3}, the multipole coefficients $M_m^{(j)}$ are obtained as follows:
$$
M_m^{(1)} = (-1)^{m+1} M_m^{(2)}=\frac{R^{m+1}}{2m}\sum_{n\neq 0}\frac{ \alpha |n| e^{-|n|s} }{\lambda- \lambda_n^+} b_{nm}.
$$

Suppose the distance $d$ is large compared to the radius $R$. Let us set $R=1/2$ for simplicity. Then the parameters have the following asymptotic behaviors:
$$
\alpha \sim d, \quad s \sim \ln d. 
$$
In this regime, one can  see that
$$
|M_m^{(1)}|=|M_m^{(2)}| \lesssim d^{-(m-1)},
$$
which shows that, in the limit $d\rightarrow \infty$, only the dipole term ($m=1$) remains non-zero.  So, when the distance between the disks is large, each disk can be approximated by a dipole source. 
However, as the disks get closer, the multipole coefficients decay very slowly as $m$ increases. The interaction between the disks become stronger and the higher order multipoles play important roles. Moreover, in this case, each of the disks cannot anymore be approximated by a dipole.

In the close-to-touching case, the solution by the spectral decomposition given in section \ref{subsec:exact} is much more efficient. It converges much faster than the multipole expansions. However, we remark that the hybridization scheme can be easily extended to multi particles system with more than 3 disks and it is highly efficient when the disks are well separated.

\section{Multipole expansion with respect to origin } \label{sec_far}

In this section, we compute the far field expansion of the scattered field $u-H$ for two disks. Specifically, we calculate the Contracted Generalized Polarization Tensors (CGPTs), which are building blocks of the expansion of the scattered field.

We will compare the first order CGPTs of two disks with that of an ellipse. It is known that, when the permittivity of the particle is positive, the GPTs of two disks can be well approximated by that of a carefully chosen ellipse, which is called the equivalent ellipse. However, we shall see that in the case of plasmonic particle, whose permittivity is negative, this does not hold any longer. Therefore, the spectral theory of the Neumann-Poincar\'{e} operator is essential in investigating the plasmonic interaction between the multi particles.
 
\subsection{CGPTs and Multipole expansion}

For a positive integer
$m$, let $P_m(x)$ be the complex-valued polynomial
\begin{equation}
P_m(x) = (x_1 + ix_2)^m := r^m \cos m\theta + i r^m \sin m\theta. \label{eq:Pdef}
\end{equation}
where $(r,\theta)$ are  the  polar coordinates.
We introduce the {\it generalized polarization tensors} by \cite{book3}
\begin{align*}
M^{cc}_{mn} = \int_{\partial \Omega} \Re \{ P_n\}  (\lambda I - \mathcal{K}_\Omega^*)^{-1} [\frac{\partial \Re \{ P_m\}}{\partial \nu}]\, d\sigma,  \\
M^{cs}_{mn} = \int_{\partial \Omega} \Im \{ P_n\}  (\lambda I - \mathcal{K}_\Omega^*)^{-1} [\frac{\partial \Re \{ P_m\}}{\partial \nu}]\, d\sigma,\\
M^{sc}_{mn} = \int_{\partial \Omega} \Re \{ P_n\}  (\lambda I - \mathcal{K}_\Omega^*)^{-1} [\frac{\partial \Im \{ P_m\}}{\partial \nu}]\, d\sigma,\\
M^{ss}_{mn} = \int_{\partial \Omega} \Im \{ P_n\}  (\lambda I - \mathcal{K}_\Omega^*)^{-1} [\frac{\partial \Im \{ P_m\}}{\partial \nu}]\, d\sigma. 
\end{align*}

Since $H$ is harmonic in $\mathbb{R}^2$, it can be written as
\beq
    H (x) = \mbox{(const.)}+\sum_{n=1}^\infty h_n^c r^n \cos n \theta + h_n^s \sin n \theta.
\label{series}
\eeq
An addition formula for the fundamental solution $\Gamma$ yields \cite{book3}
 \beq\label{eq-green GPT expansion2}
\Gamma(x-y) = \sum_{m=0}^{\infty} \frac{(-1)}{2\pi m}  \frac{\cos(m \theta)}{ r^m} r_{{y}}^m\cos(m \theta_{{y}}) + \frac{(-1)}{2\pi m} \frac{\sin(m \theta)}{ r^m} r_{{y}}^m\sin(m \theta_{{y}})
\eeq
for $|x|>|y|$.
Here $(r_y,\theta_y)$ are the polar coordinates of $y$, i.e., $y=r_y e^{i\theta_y}$.

We then have from \eqref{scattered} and the definitions of the CGPTs that
\begin{align}
(u-H)(x) &= \sum_{m=1}^\infty  \frac{(-1)}{2\pi m}  \frac{\cos(m \theta)}{ r^m} (M_{nm}^{cc} h_n^c + M_{nm}^{sc} h_n^s)\nonumber
\\
&\quad
+ \sum_{m=1}^\infty\frac{(-1)}{2\pi m}  \frac{\sin(m \theta)}{ r^m} (M_{nm}^{sc} h_n^c + M_{nm}^{ss} h_n^s)  \label{EF2}
\end{align}
for large enough $|x|$. Note that the above expression is a multipole expansion with respect to the origin $(0,0)$. Once we get the CGPTs explicitly, the above formula can give highly efficient and accurate approximations for the scattered field when $|x|$ is large.


\subsection{Explicit computations of CGPTs for two disks}
Here we explicitly compute the CGPTs associated with two disks.
It is easy to see that the CGPTs associated with $B_1\cup B_2$ can be represented as
\begin{align*}
M^{cc}_{nm} = \big\langle \widetilde{P}_m^c , (\lambda I-\mathbb{K}^*)^{-1}[\p_\nu \widetilde{P}_n^c] \big\rangle_{L^2}.
\\
M^{cs}_{nm} = \big\langle \widetilde{P}_m^s , (\lambda I-\mathbb{K}^*)^{-1}[\p_\nu \widetilde{P}_n^c] \big\rangle_{L^2}.
\\
M^{sc}_{nm} = \big\langle \widetilde{P}_m^c , (\lambda I-\mathbb{K}^*)^{-1}[\p_\nu \widetilde{P}_n^s] \big\rangle_{L^2}.
\\
M^{ss}_{nm} = \big\langle \widetilde{P}_m^s , (\lambda I-\mathbb{K}^*)^{-1}[\p_\nu \widetilde{P}_n^c] \big\rangle_{L^2},
\end{align*}
where
\begin{align*}
\widetilde{P}_m^c = \begin{bmatrix}
 \Re\{P_m\} |_{\p B_1} \\
   \Re\{P_m\} |_{\p B_2}
 \end{bmatrix} =  \begin{bmatrix}
r^m \cos m \theta  |_{\p B_1} \\
r^m \cos m \theta|_{\p B_2}
 \end{bmatrix},
\end{align*}
and
\begin{align*}
\widetilde{P}_m^s = \begin{bmatrix}
 \Im\{P_m\} |_{\p B_1} \\
   \Im\{P_m\} |_{\p B_2}
 \end{bmatrix}=  \begin{bmatrix}
r^m \sin m \theta  |_{\p B_1} \\
r^m \sin m \theta|_{\p B_2}
 \end{bmatrix}.
\end{align*}
Due to the symmetry of $B_1 \cup B_2$, we have
$$
M^{cs}_{nm}=M^{sc}_{nm}=0.
$$

Let us compute $M_{nm}^{cc}$ and $M_{nm}^{ss}$. By the spectral decomposition of $\mathbb{K}^*$, we have
\begin{align*}
M_{nm}^{cc}&= \sum_{k \neq 0} \frac{\la \widetilde{P}_m^c, \Psi_{k}^+\ra_{L^2}\la \Psi_{k}^+,\p_\nu\widetilde{P}_n^c \ra_*}{\lambda-\lambda_{k}^+}
+\sum_{k \neq 0} \frac{\la \widetilde{P}_m^c,\Psi_{n}^-\ra_{L^2}\la \Psi_{k}^-,\p_\nu\widetilde{P}_n^c \ra_*}{\lambda-\lambda_{k}^-}
\\
&=\sum_{k \neq 0} \frac{\left(\frac{1}{2}-\lambda_{k}^+\right)\la \widetilde{P}_m^c,\Psi_{k}^+\ra_{L^2} 
\overline{\la \widetilde{P}_n^c,\Psi_{k}^+\ra_{L^2}}
}{\lambda-\lambda_{k}^+}
+\sum_{k \neq 0} \frac{\left(\frac{1}{2}-\lambda_{k}^-\right)\la \widetilde{P}_m^c,\Psi_{k}^-\ra_{L^2} 
\overline{\la \widetilde{P}_n^c,\Psi_{k}^-\ra_{L^2}}
}{\lambda-\lambda_{k}^-}.
\end{align*}
So we need to compute $\la \widetilde{P}_n^c,\Psi_{k}^\pm \ra_{L^2}$.
Straightforward but tedious computations show that
\begin{align*}
\la \widetilde{P}_m^c,\Psi_{k}^+\ra_{L^2} &= \frac{\pi}{2} \alpha^m \sqrt{\frac{|k|}{\frac{1}{2}-\lambda_k^+ }} e^{-|k|s} \widetilde F(m,|k|),
\end{align*}
where $\widetilde F(m, |k|)$ is defined by 
\begin{equation} \label{defF}
\widetilde F(m,k)=\sum_{l=0}^{\min(k,m)} (-1)^m \begin{pmatrix}
m \\ l
\end{pmatrix}
\begin{pmatrix}
k-l+m-1
\\
k-m
\end{pmatrix}.
\end{equation}

By a symmetry consideration, we have
$$
\la \widetilde{P}_m^c,\Psi_{k}^-\ra_{L^2} = 0.
$$

\begin{theorem}
The following explicit formulas for the CGPTs hold:
\begin{align*}
\ds M_{nm} 
&=\begin{pmatrix}
M_{nm}^{cc} & M_{nm}^{cs}
\\
M_{nm}^{sc} & M_{nm}^{ss}
\end{pmatrix}
\\
  &= \ds\begin{pmatrix}\ds
 \pi a^{m+n} \sum_{k=1}^\infty \frac{  \widetilde F(m,k) \widetilde F(n,k) k e^{-2ks}  }{\lambda-\lambda_{k}^+} &  0 \\
0 &  \ds\pi a^{m+n} \sum_{k=1}^\infty \frac{  \widetilde F(m,k) \widetilde F(n,k) k e^{-2ks}  }{\lambda-\lambda_{k}^-}
 \end{pmatrix}.
\end{align*}
\end{theorem}

\begin{remark}
 The CGPTs $M_{nm}$ contain an infinite number of poles with respect to $\lambda$ and the poles are the eigenvalues of the Neumann-Poincar\'{e} operator (or the plasmon resonances). When $\lambda$ is close to one of the eigenvalues $\lambda_k^\pm$, the scattered field $u-H$ is greatly enhanced.
\end{remark}

\begin{remark}
It is known that, if $D$ is an ellipse of the form $R_\theta(B')$ where $R$ is a rotation y $\theta$ and $B'$ is an ellipse of the form $(x_1/a_1)^2 + (x_2/a_2)^2 <1 $, then the first order polarization tensor $M_{11}(\lambda,D)$ is
$$
M(\lambda,D) =R 
\begin{pmatrix}
\ds\frac{|B|}{\lambda - \frac{1}{2}\frac{a_1-a_2}{a_1+a_2}}
& 0
\\
0 &\ds \frac{|B|}{\lambda - \frac{1}{2}\frac{a_1-a_2}{a_1+a_2}}
\end{pmatrix}
R^T.
$$
Contrary to the two disks case, $M_{11}(\lambda,D)$ contains only two poles with respect to $\lambda$. So the ellipse cannot be used to get approximation of the scattered field generated by two disk when the particles are plasmonic.
\end{remark}

\section{Numerical illustrations} \label{sec6}

Here we compare the numerical results obtained from the trasnformation based solution and the plasmon hybridization (or multipole expansion) based solution. We compute the field at the origin $\nabla (u-H)(0,0)$ as a function of frequency.

We suppose that the radius $R$ is $R=1$ and the distance $d$ is $d=0.1$. Since $d/R$ is small, we expect that interaction between two disks is strong. For the permittivity $\ep_m$ of two disks, we use \eqnref{drude} with $\omega_p = 3$ and $\gamma=0.02$.

According to 
the plasmon hybridization (or multipole expansion) solution \eqnref{EF3}, the field is given by
$$
\nabla (u-H)(0,0) = (-2)\sum_{m=1}^\infty \frac{m M_m^{(1)}}{(R+d/2)^{m+1}} \mathbf{e}_1 \approx (-2)\sum_{m=M}^\infty \frac{m M_m^{(1)}}{(R+d/2)^{m+1}} \mathbf{e}_1,
$$
and, by using the transformation based solution \eqnref{u_H_spectral}, the field at the origin can be computed as
$$
\nabla (u-H)(0,0) = \frac{2}{\alpha}\sum_{n=1}^\infty \frac{2 \alpha n e^{-ns}}{\lambda-\lambda_n^+} (-1)^n \mathbf{e}_1 \approx  \frac{2}{\alpha}\sum_{n=N}^\infty \frac{2 \alpha n e^{-ns}}{\lambda-\lambda_n^+} (-1)^n \mathbf{e}_1,
$$
where $M$ and $N$ are the truncation numbers. When $N$ or $M$ increases, the accuracy will improve.

\begin{figure*}
\begin{center}

\epsfig{figure=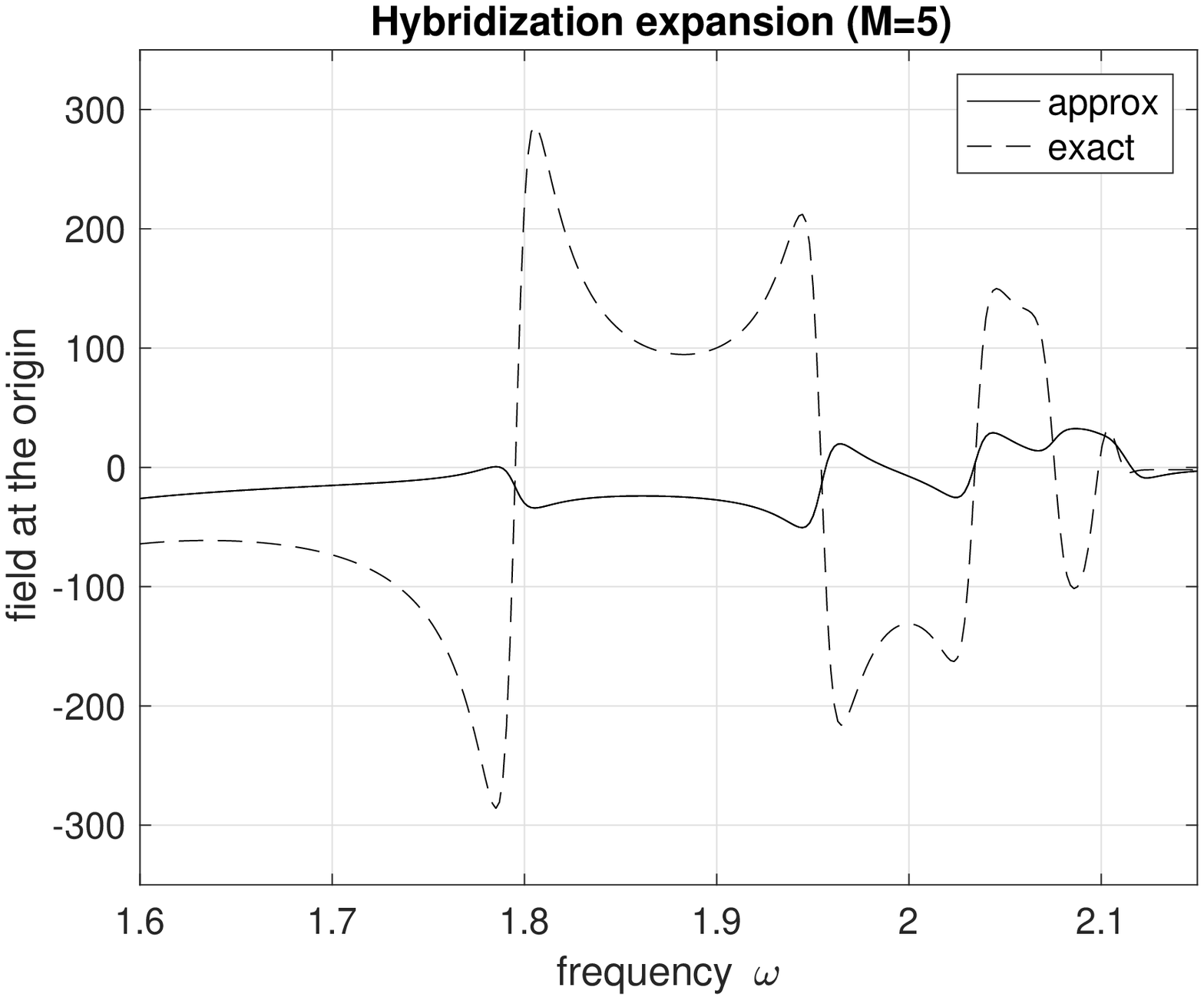,width=5cm}\hskip-0.1cm
\epsfig{figure=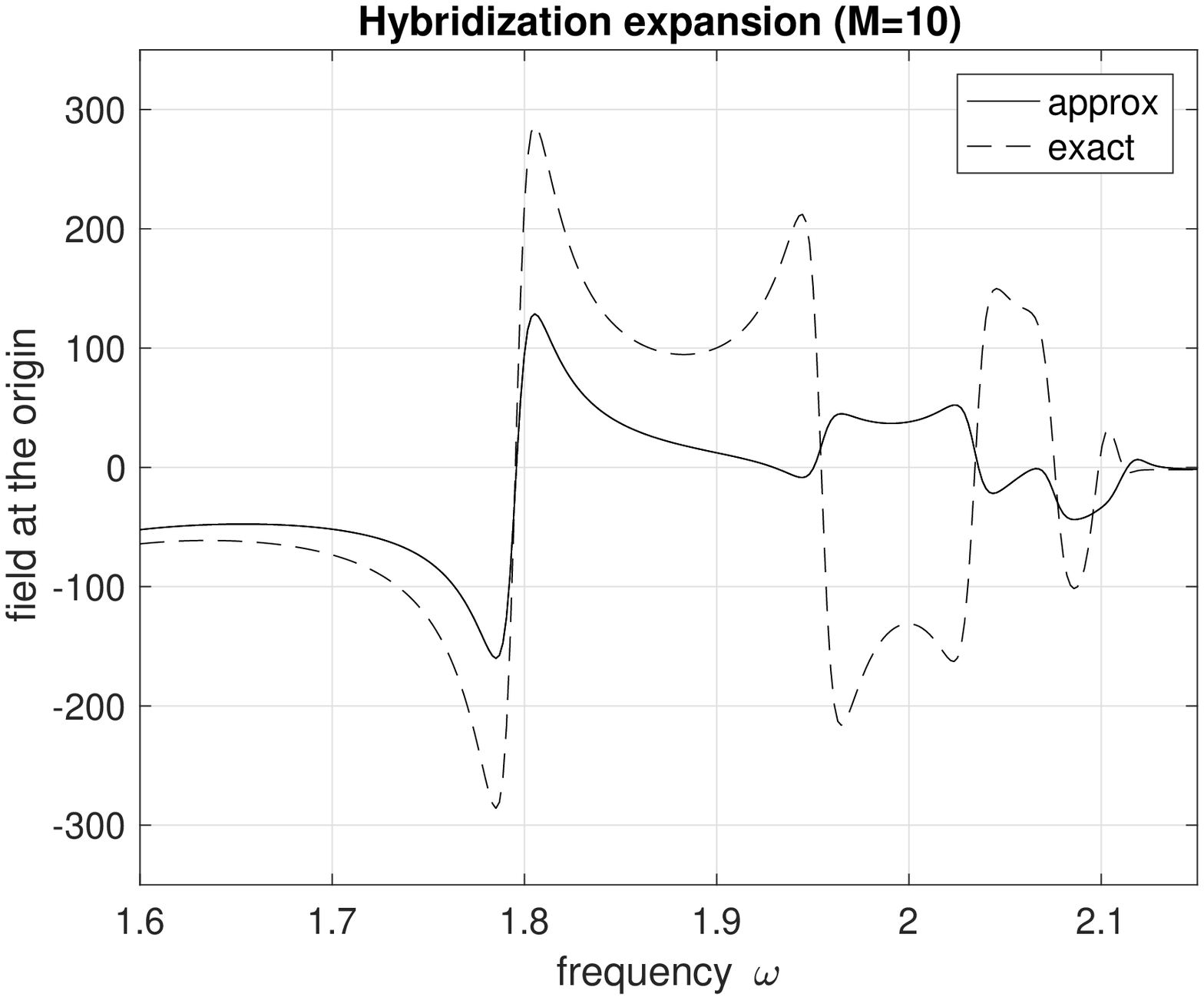,width=5cm}\hskip-0.1cm
\epsfig{figure=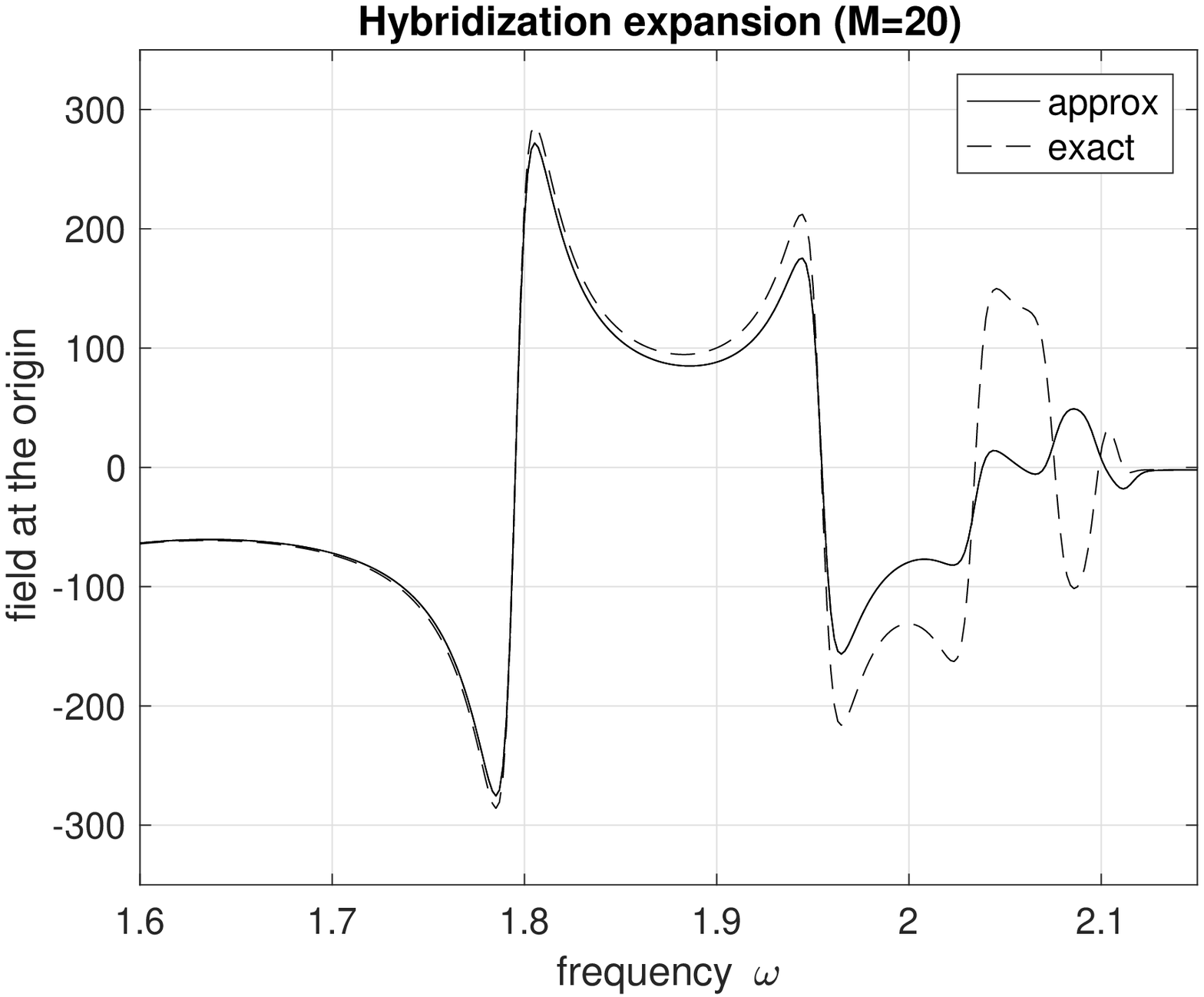,width=5cm}

\vskip.5cm

\epsfig{figure=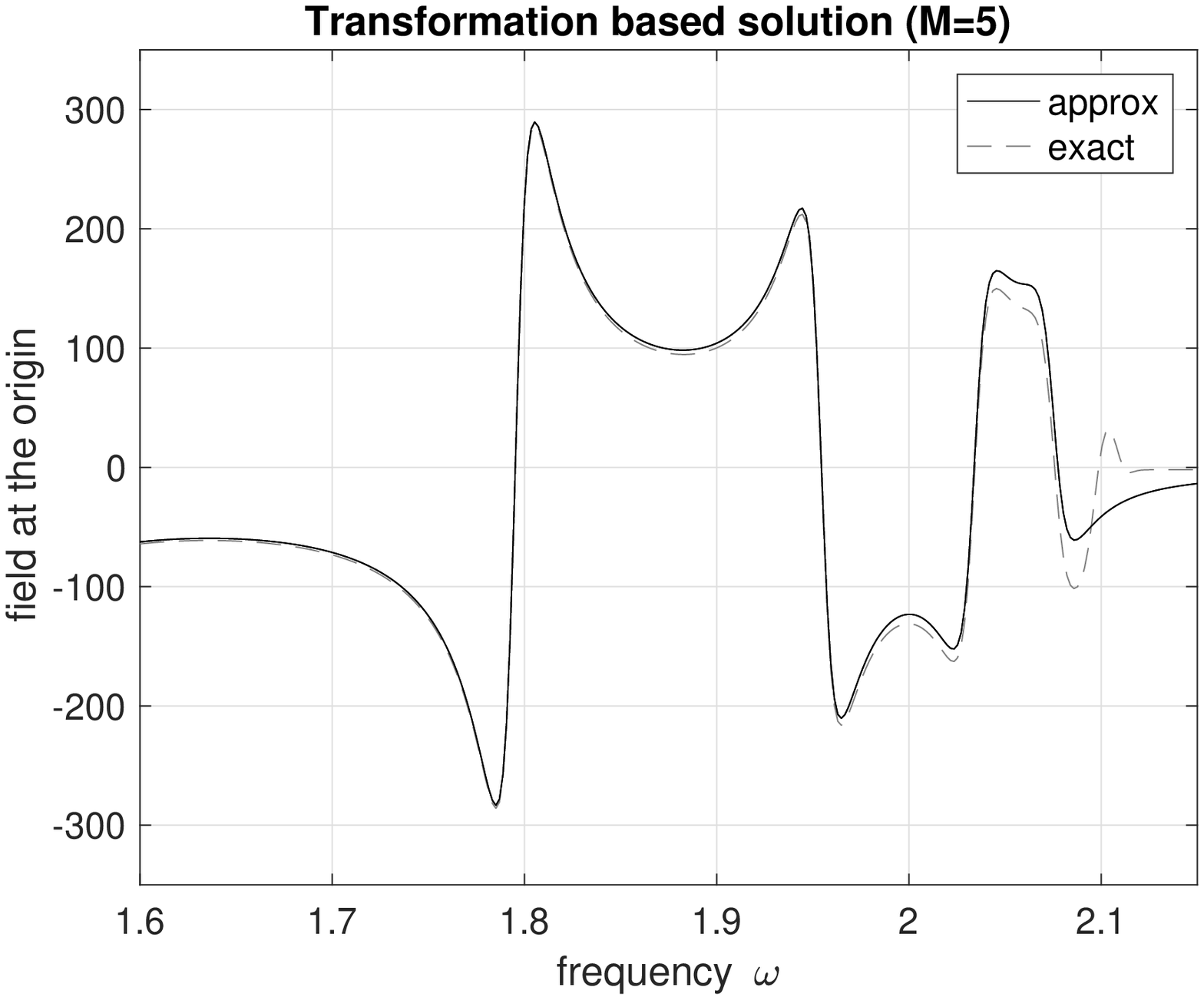,width=5cm}\hskip-0.1cm
\epsfig{figure=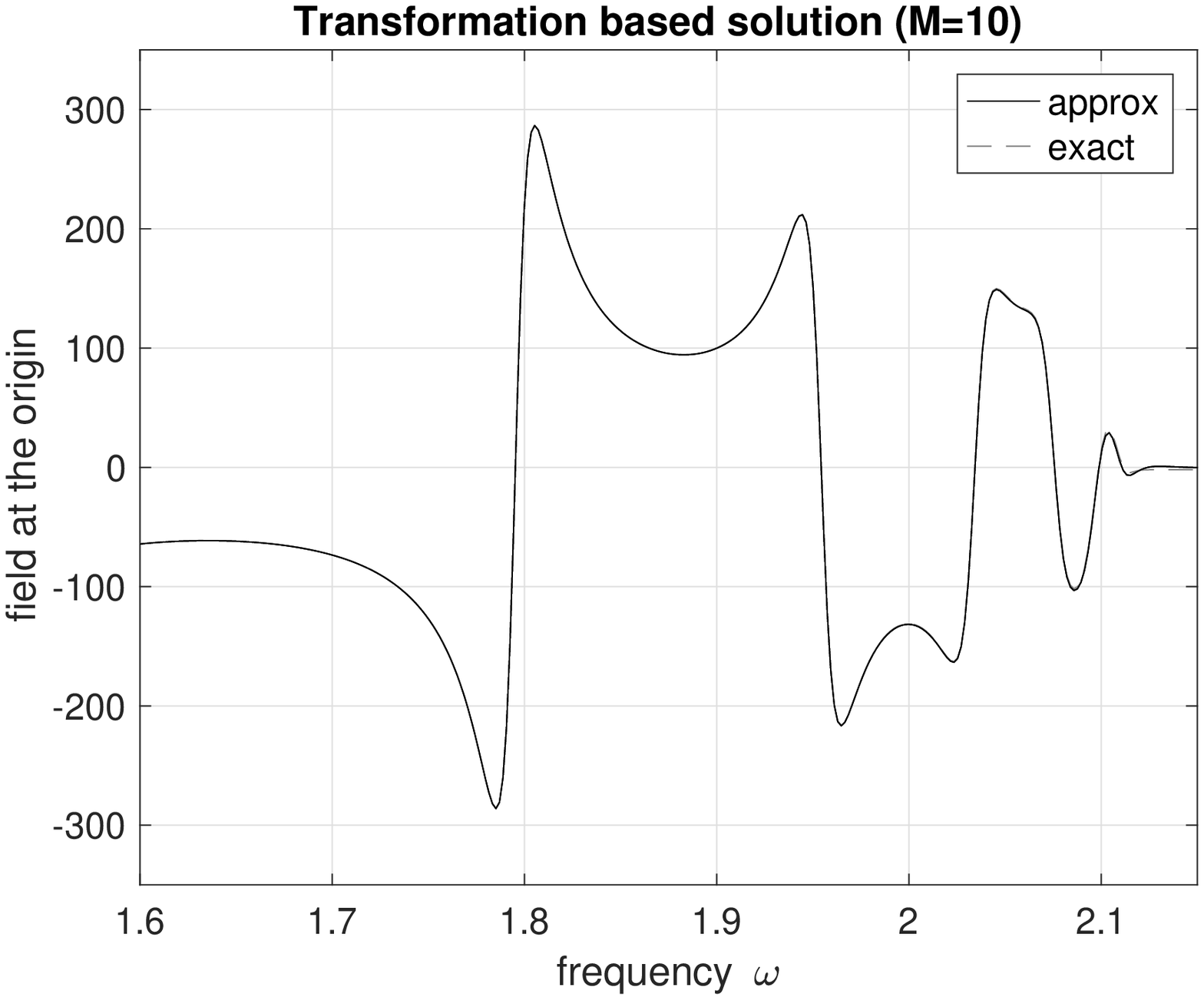,width=5cm}\hskip-0.1cm
\epsfig{figure=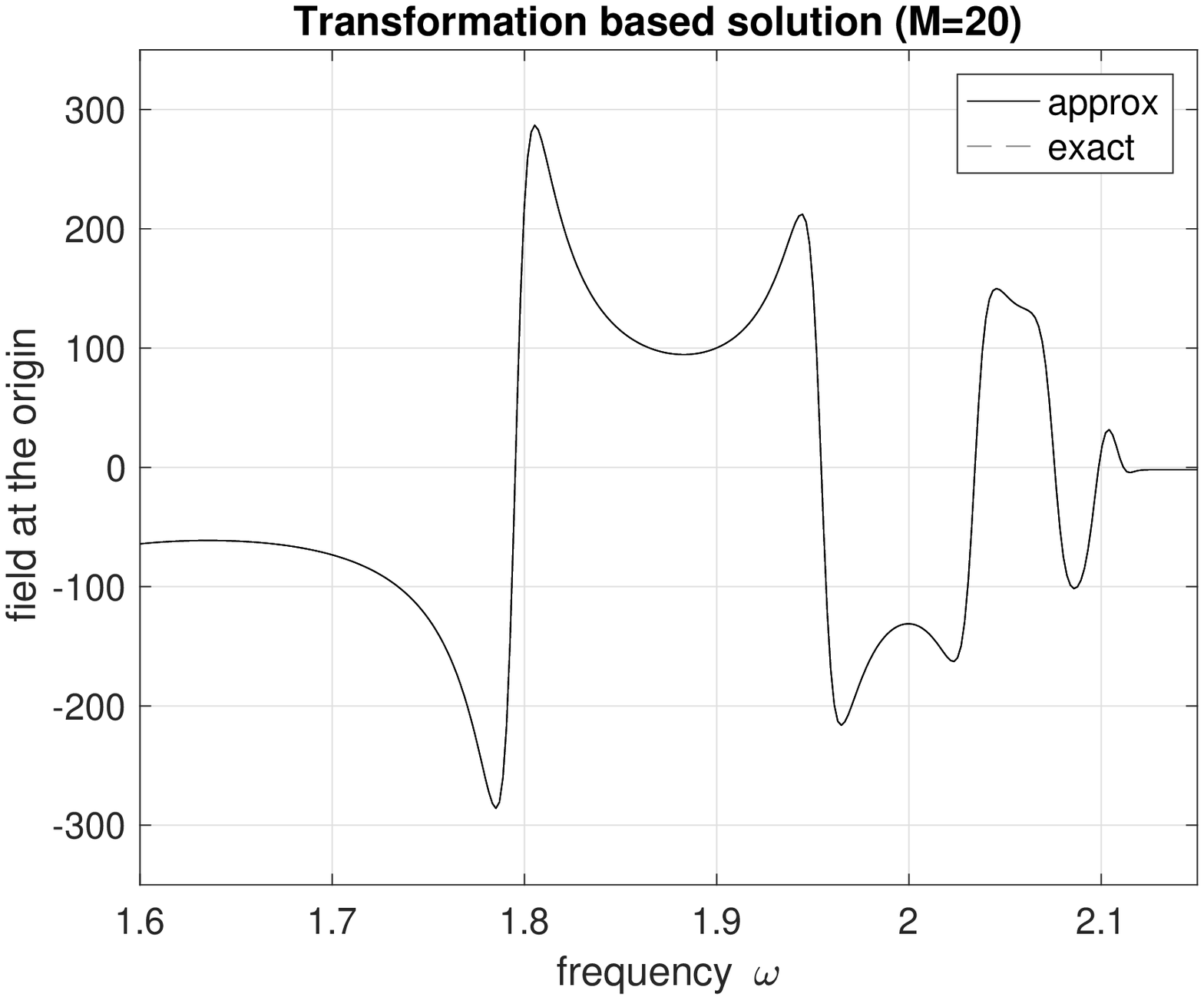,width=5cm}
\end{center}
\caption{The field at the origin as a function of frequency. The multipole expansion method for $M=5,10$ and $20$ (first row), and the transformation based solution for $N=5,10$ and $20$ (second row).}
\label{fig:field}
\end{figure*}

In Figure \ref{fig:field}, we compute $\mathbf{e}_1\cdot\nabla (u-H)(0,0)$ as a function of frequency. In the first row (or the second row), we show the numerical results computed by the plasmon 
hybrdized expansion (or the transformation based solution) when $M=5,10,20$ (or $N=5,10,20$), respectively. 
It indicates that the transformation based method is much more efficient. The multipole expansion method gives an inaccurate result even when we use high number of orders with $M=20$. On the contrary, the accuracy of the transformation based method is pretty good although we use a low order with $N=5$.


\begin{thebibliography}{99}

\bibitem{alu}
A. Alu and N. Engheta, 
Dynamical theory of artificial optical magnetism produced by rings of plasmonic nanoparticles, 
Physical Review B, 78, 085112(2008).


\bibitem{ACKLM2} H. Ammari, G. Ciraolo, H. Kang, H. Lee, and G.W. Milton, Spectral theory
of a Neumann-Poincar\'e-type operator and analysis of anomalous localized resonance II, Contemp. Math., 615 (2014), 1--14.


\bibitem{pierre}
 H. Ammari, Y. Deng, and P. Millien,
Surface plasmon resonance of nanoparticles and applications in imaging, Arch. Ration. Mech. Anal., 220 (2016), 109--153. 


\bibitem{matias} H. Ammari, P. Millien, M. Ruiz, and H. Zhang, 
Mathematical analysis of plasmonic nanoparticles: the scalar case,  Archive on Rational Mechanics and Analysis, 224 (2017), 597--658.


\bibitem{sinum}  H. Ammari, H. Kang, E. Kim, and M. Lim, Reconstruction of closely spaced small inclusions, SIAM J. Numer. Anal., 42 (2005),  2408--2428.


\bibitem{hh}
H. Ammari, M. Ruiz, S. Yu, and H. Zhang, 
Reconstructing fine details of small objects by using plasmonic spectroscopic data. To appear in SIAM Journal on Imaging Sciences.

\bibitem{book3}  H. Ammari, J. Garnier, W. Jing, H. Kang, M. Lim, K. S\o lna, and H. Wang,
\textsl{Mathematical and Statistical Methods for Multistatic Imaging}, Lecture Notes in Mathematics, Volume 2098, Springer, Cham, 2013.




%
%

\bibitem{mihai} H. Ammari, M. Putinar, M. Ruiz, S. Yu, and H. Zhang, Shape reconstruction of nanoparticles from their associated  plasmonic resonances, J. Math. Pures Appl., DOI:10.1016/j.matpur.2017.09.003. 


\bibitem{matias2} H. Ammari, M. Ruiz, S. Yu, and H. Zhang,
Mathematical analysis of plasmonic resonances for nanoparticles: the full Maxwell equations, 
Journal of Differential Equations, 261 (2016), 3615--3669. 



%
%
%
\bibitem{kang1} K. Ando and H. Kang, Analysis of plasmon resonance on smooth domains using spectral properties of the Neumann-Poincar\'e operator, J. Math. Anal. Appl., 435 (2016), 162--178.


\bibitem{anker}
J. N. Anker, W. P. Hall, 
O. Lyandres, N. C. Shah, J. Zhao, and 
R. P. Van Duyne, 
Biosensing with plasmonic nanosensors, Nature material, 7 (2008), 442--453.


%


%
%
\bibitem{BT_disks}
E. Bonnetier and F. Triki,
Pointwise bounds on the gradient and the spectrum of the Neumann–Poincar\'{e} operator: the case of 2 discs, Contemp. Math 577, 81--92.


\bibitem{triki} E. Bonnetier and F. Triki, On the spectrum of the Poincar\'e variational problem for two close-to-touching inclusions in 2D, Arch. Ration. Mech. Anal.,  209  (2013),   541--567.


\bibitem{dmitri2013}
D.K. Gramotnev1 and S.I. Bozhevolnyi, 
Nanofocusing of electromagnetic radiation, Nature Photonics,  8 (2014), doi: 10.1038/NPHOTON.2013.232.
%
\bibitem{Gri12}
D. Grieser, The plasmonic eigenvalue problem, Rev. Math. Phys.,  26  (2014),   1450005.

\bibitem{hybrid} N.J. Halas, S. Lal, W.S. Chang, S. Link, and P. Nordlander, Plasmons in strongly coupled metallic nanostructures, Chemical Rev., 111 (2011), 3913--3961.   

%
\bibitem{plasmon4} P.K. Jain, K.S. Lee, I.H. El-Sayed, and M.A. El-Sayed, Calculated absorption and scattering properties of gold nanoparticles of different size, shape, and composition: Applications in biomedical imaging and biomedicine, J. Phys. Chem. B, 110 (2006), 7238--7248.
%


%
%
%
%
%
%
%
%
%
\bibitem{kelly} K.L. Kelly, E. Coronado, L.L. Zhao, and G.C. Schatz, 
The optical properties of metal nanoparticles: The influence of size, shape, and dielectric environment, J. Phys. Chem. B, 107 (2003), 668--677. 

\bibitem{lauchner}
A. Lauchner, A.E. Schlather, A. Manjavacas, Y. Cui, M.J. McClain, G.J. Stec, F.J. Garcia de Abajo, P. Nordlander, and N.J. Halas, 
Molecular Plasmonics,
Nano Letters 2015 15 (9), 6208--6214
DOI: 10.1021/acs.nanolett.5b02549.

\bibitem{boris2010}
B. Luk'yanchuk, N.I. Zheludev,	S.A. Maier, N.J. Halas,	P. Nordlander, H. Giessen, and C.T. Chong,
The Fano resonance in plasmonic nanostructures and metamaterials,
Nature Materials 9, 707--715 (2010) doi:10.1038/nmat2810.




%
%
%
%
\bibitem{plasmon1} I.D. Mayergoyz, D.R. Fredkin, and Z. Zhang, Electrostatic (plasmon) resonances in nanoparticles, Phys. Rev. B, 72 (2005), 155412.
%



%
\bibitem{SC10}
D. Sarid and W.A. Challener, \textsl{Modern Introduction to Surface Plasmons: Theory, Mathematical Modeling, and Applications},
Cambridge University Press, New York, 2010.


\bibitem{simovski}
C.R. Simovski, A.J. Viitanen, and
S.A. Tretyakov, 
Sub-wavelength resolution in linear arrays of plasmonic particles,
Journal of Applied Physics 101, 123102 (2007), doi: http://dx.doi.org/10.1063/1.2745315.


\bibitem{yu17}
S. Yu and H. Ammari, 
Plasmonic interaction between nanospheres. To appear in SIAM Review.


\end{thebibliography}
\end{document}